%% file: two_theorems_2405.tex
\documentclass[12pt,reqno]{article}
\usepackage{amsmath,amsfonts,amssymb,amsthm,indentfirst}

\input
cd.tex 
\newtheorem{thm}{Theorem} \newtheorem{prob}{Problem}

\newtheorem{prop}{Proposition}[section]

\theoremstyle{definition} \newtheorem{defn}{Definition}[section]
\newtheorem*{rem*}{Remark} \newtheorem{remark}{Remark}

\def\GL{\textrm{GL}}

           \def\int{\textrm{int}}

\def\Ker{\textrm{Ker}}

\def\HP{\textrm{HP}}
\def\GL{\textrm{GL}}
\numberwithin{equation}{section}
\date{}

\title{Notes on Engel groups and Engel  elements in groups. Some generalizations}
\author{Boris Plotkin, Hebrew University, Jerusalem}

\begin{document}

\maketitle
\begin{abstract}
Engel groups and Engel elements became popular in 50s. We consider
in the paper the more general nil-groups and nil-elements in
groups. All these notions are related to nilpotent groups and
nilpotent radicals in groups. These notions generate problems
which are parallel to Burnside problems for periodic groups.

The first three theorems of the paper are devoted to nil-groups
and Engel groups, while the other results are connected with the
further generalizations. These generalizations extend the theory
to solvable groups and solvable radicals in groups. The paper has
two parts. The first one (sections 2-4) deals with old ideas,
while the second one (sections 5-9) is devoted to generalizations.

\end{abstract}

\section{Introduction}

 This paper is a nostalgic  reminiscence on group theory
of 50s (just last century). In some sense this feedback to the
past is inspired by the paper \cite{BGGKPP} and by the recent
talks on $PI$-algebras by L. Rowen and A. Kanel-Belov.
 Recall some
definitions and some necessary old results.

We distinguish Engel groups and nil-groups, Engel elements and
nil-elements \cite{Pl3}.

Let $F_2=F(x,y)$ be the free group. Define the sequence $$
e_1(x,y)=[x,y], e_2(x,y)=[e_1(x,y),y],\ldots,
e_n(x,y)=[e_{n-1}(x,y),y], $$ where $x,y\in F_2$. Let now $G$ be
an arbitrary group.
 \begin{defn} An element $g\in
G$ is called nil-element if for every $a\in G$ there is $n=n(a,g)$
such that $e_n(a,g)=1$.
\end{defn}
\begin{defn}
A group $G$ is called {\it nil-group} if every its element is
nil-element.
\end{defn}
Every locally-nilpotent group is a nil-group, but the opposite
 is not
true (\cite{GS}).
\begin{defn}
A group $G$ is called {\it Engel group} if it satisfies an
identity $e_n(x,y)\equiv 1$ for sone $n$.
\end{defn}

In this case we call the group $n$-Engel. The variety of $n$-Engel
groups is denoted by $E_n$ and let $F^n_k$ be the free group with
$k$ free generators in this variety .

There is a long-standing conjecture that the group $F^n_k$ is not
nilpotent but up to now there are no reasonable approaches to
this problem. We show (Theorem 1) the restricted solution of this
problem similar to a solution of the restricted Burnside problem
(see for example \cite{Ne},\cite{Ko}, etc.)
 \begin{defn} An
element $g\in G$ is called Engel element, if there exists
$n=n(g)$ such that for every $x\in G$ the identity
$e_n(x,g)\equiv 1$ holds in $G$.
\end{defn}

Thus, the definition of an Engel element differs from the
definition of a nil-element. However, following the tradition
sometimes we use the term Engel element also for nil-elements
(with the meaning unbounded Engel elements).

The conditions on a group $G$ provided the set of all nil-elements
or (and)  the set of all Engel elements constitute a subgroup are
considered. In view of the latter problem note the following
general result \cite{Pl2}:

Let the group $G$ has an ascending normal series with the locally
neotherian quotients. Then the set of all its nil-elements
constitute the subgroup in $G$ coinciding with the locally
nilpotent radical  $\HP(G)$.

This Theorem has been preceded by the  similar theorem for the
case when the factors of the normal series are locally nilpotent
\cite {Pl1} and the theorem of Baer \cite{Ba}, stating that the
nilpotent radical of a Noetherian group coincides with the set of
its nil-elements.

Theorem of Baer follows from the Lemma \cite{Pl2} we are going to
recall:

Let $G$ be an arbitrary group, $g$ its nil-element. Then there
exists in $G$ a normal series of nilpotenet subgroups

$$
H_1\subset H_2\ldots, H_n\subset\ldots,
$$
 where $H_1=\{g_1\}$, $H_n=\{H_{n-1},h_ngh_n^{-1}\}$ for some
 $h_n\in G$. Here and elsewhere $\{\ \}$ stands for the subgroup
 generated by some elements. This series stops on some place $n$
 if $H_n$ is a normal subgroup in $G$.

Note here the following result of A.~Tokarenko \cite{To1}, see
also \cite{Pl3}:

Let a group $G$ be a subgroup in some $\GL_n(K)$, where $K$ is a
commutative ring with 1. Then the set of all nil-elements in $G$
is the locally nilpotent radical $\HP(G)$.

Now we formulate two  theorems of this paper:

\begin{thm} In any variety $E_n$ all its locally nilpotent groups
form a subvariety.
\end{thm}

In the second theorem we consider $PI$-groups that is the groups
which can be embedded to a group of invertible elements of some
$PI$- algebra over a filed $P$.

For example, the full matrix group $\GL_n(P)$ and all its
subgroups are $PI$-groups.

It has been proved by Procesi \cite{Pr} and Tokarenko \cite{To2}
that every periodic $PI$-group is locally finite. The following
theorem has the similar flavor:

\begin{thm} Every nil-$PI$-group $G$ is locally nilpotent.
\end{thm}

In fact we will prove the more general result:

\begin{thm} In every $PI$-group $G$ the set of all its nil-elements coincides
with the locally nilpotent radical $\HP(G)$.
\end{thm}

Since other results will require additional definitions, they will
be formulated later.

Engel groups and Engel elemets in groups are related to nilpotent
groups and nilpotent radicals in groups. Along with these elements
we consider also their Engel-like generalizations, which are tied
with solvable groups and solvable radicals in  groups.


 Note also two facts which will be  used in the sequel.

 First of all this is the theorem by Wilson (\cite{Wi}) which
 states that every residually finite finitely generated Engel
 group is nilpotenet.

 Second, we use the following Kaluzhnin's theorem \cite{Ka}. Let the group
 $G$ acts unitriangularly  and faithfully in the space $V$. Then the
 group $G$ is nilpotent. Unitriangularity means that there is a
 series
 $$
 V=V_0\supset V_1\supset\ldots\supset V_n=0
 $$
 in $V$ such that all members of the series are invariant in respect to the action of the
 group $G$ and in all quotients of the series the group $G$ acts
 trivially. If the group $G$ acts in $V$ faithfully then $G$ has
 the nilpotency class $n-1$.

\section{Proof of the theorem 1. }

Prove first the following:

\begin{prop} In the variety $E_n$ for any natural $k$ there exists
a nilpotent group with $k$ generators $\widetilde F^n_k$ such that
every $k$-generated nilpotent group $G\in E_n$ is a homomorphic
image of the group $\widetilde F^n_k$.
\end{prop}

Proof. Let us start with the free in $E_n$ group $F^n_k$. Let
$H^n_k$ be the intersection of all $H\lhd F^n_k$ with nilpotent
$F^n_k/H$. The group $\widetilde F^n_k=F^n_k/H^n_k$ is residually
nilpotent. This group is also residually finite since every
finitely generated nilpotent group is residually finite. Besides,
this group is Engel. By the result from \cite{Wi} the group
$\widetilde F^n_k$ is nilpotent.

Let $G$ be an arbitrary group in $E_n$ with $k$ generators. There
is a surjection $F^n_k\to G$. Let $H$ be the kernel of this
homomorphism. Then $H\supset H^n_k$ and this gives a surjection
$\widetilde F^n_k\to G.$ $\square$

This proposition is equivalent, in fact, to the theorem 1. Indeed,
let $\Theta$ be the class of all locally nilpotent groups in
$E_n$. This class is closed in respect to taking subgroups and
homomorphic images. Let us check that the class is closed in
respect to Cartesian products.

Let $G=\prod_{\alpha\in I}G_\alpha$ where all $G_\alpha$ are
locally nilpotent groups in the variety $E_n$. Take a finitely
generated subgroup $H$ in $G$ with $k$ generators. The group $H$
is approximated by $k$ generated nilpotent groups $H_\alpha\subset
G_\alpha$. The nilpotency class of all $H_\alpha$ is bounded by
the nilpotency class of the group $\widetilde F^n_k$. Thus, $H$ is
nilpotent and $G$ is locally nilpotent.

This means that the class $\Theta$ is a variety. It is easy to see
that the group $\widetilde F^n_k$ is the free group with $k$
generators in the variety $\Theta$.

\section{$PI$-groups}


Let us fix a field $P$ and for every group $G$ consider a
representation $G\to A,$  where $A$ is an associative algebra
with 1 over the ground field $P$ and the arrow means a
homomorphism of the group $G$ to the group of invertible elements
of the algebra $A$. If this representation is faithful then we
say that the algebra $A$ is a linear envelope of the group $G$.
The group algebra $PG$ is the universal linear envelope. We
consider  the groups $G$ from the point of view of the possible
linear envelopes. In particular, $G$ is a linear group if it has
a finite dimensional linear envelope.

 \begin{defn} A group $G$ is a $PI$-group if it has a linear envelope $A$ which
 is a $PI$-algebra.
 \end{defn}

 Let us fix $A$ and let $G_0$ be the group of invertible elements of the algebra
 $A$. We consider $G$ as a subgroup in $G_0$.

 In a $PI$-algebra $A$ consider a series of ideals:

 $$
 U_0=0\subset U_1\subset U_2\subset A,
 $$
 where $U_1$ is the sum of the nilpotent ideals of $A$ and  $U_2$ is the
 Levitzky radical of $A$. It is known \cite{J} that $U_2/U_1$ is the nilpotent algebra
 and there is a faithful embedding
 $A/U_2\to M_n(K)$. Here $M_n(K)$ is the matrix algebra of the dimension $n$ and $K$ is a
 commutative ring with 1 which is a cartesian sum of fields. The
 group of invertible elements of $M_n(K)$ is $\GL_n(K)$.

 \begin{prop} Let $G$ be a $PI$-group. Then there is a chain of
 normal subgroups
 $$ 1=H_0\subset H_1\subset H_2\subset G, $$
where $H_1$ is generated by the nilpotent normal subgroups in
$G$, $H_2$ is locally nilpotent, and there is a faithful
embedding $$ G/H_2\to \GL_n(K),$$ where $K=\bigoplus P_\alpha$,
$P_\alpha$ is a field. \end{prop} Proof. First recall the known
things. Let $A$ be an associative algebra with 1, $G_0$ the group
of invertible elements in $A$, $G$ a subgroup in  $G_0$. The
group $G$ acts in the space $A$ by the rule: $a\to ag$, $a\in A$,
$g\in G$.

Let $U$ be a two-sided ideal in $A$ and $\mu:A\to A/U$ the
natural homomorphism. It induces the representation $\mu:G\to
A/U$. Then $\mu_0(1)$ is the coset $1+U$ and the kernel of
$\mu_0$ is the set of elements $g\in G$, such that $g-1\in U$,
i.e., $g\in 1+U$. We have $\Ker(\mu_0)=G\bigcap (1+U)$. The group
$G$ acts also in $A/U$ with the same kernel $G\bigcap (1+U)$.

Consider the coset $1+U$ and let $U$ be a locally nilpotent ideal.
Check that $H=1+U$ is a locally nilpotent normal subgroup in
$G_0$. The set $H$ is closed under multiplication. Let $a\in U$
and $a^n=0$. We have $(1+a)(1-a+a^2+\ldots+(-1)^{n-1}a^{n-1})=1$
and $1+a$ is an invertible element. Thus, $H=1+U$ is a subgroup
in $G_0$. This subgroup is normal since it coincides with the
kernel $G_0\to A/U$ . It remains to check that $H$ is locally
nilpotent.

Consider first the case when $U$ is a nilpotent ideal. Consider a
series $$ U=U_0\supset U_1\supset\ldots U_k\supset\ldots\supset
U_n=0, $$ where $U_k$ is a linear combination of the elements of
the form $$ a(g_1-1)\ldots(g_k-1), a\in U, g_i\in H .$$ The
series above is invariant under the action of $H$ and $H$ acts
trivially in the factors. Besides, $H$ acts trivially in $A/U$.
This means that $H$ acts in $A$ unitriangularly and faithfully.
Thus, by Kaluzhnin's theorem $H$ is nilpotent.

Let now $U$ be a locally-nilpotent ideal, $H_0=\{
g_1,\ldots,g_n\}$ be a finitely generated subgroup in $H$. Assume
that for every generator $g_i$ its inverse belongs to the set $\{
g_1,\ldots,g_n\}$. Then $g_i$, $i=1,2,\ldots,n$ generates $H_0$
as a semigroup.

For every $g_i$ take $a_i=g_i-1$. Generate by the elements
$a_1,\ldots,a_m$ a subalgebra $U_0$ in $U.$ The subalgebra $U_0$
is nilpotent. It is easy to see that $g-1\in U_0$ for every $g\in
H_0$ and $H_0\subset 1+U_0$. Take a subalgebra $U^*_0=\{U_0,1\}$
in $A$. Here $U_0$ is the nilpotent ideal in $U^*_0$. The group
$1+U_0$ acts in $U^*_0$ faithfully and unitriangularly. Hence,
$1+U_0$ is nilpotent and $H_0$ is also nilpotent. Thus, $H$ is
locally nilpotent.

Return to the situation when $A$ is a $PI$-algebra and let $$
U_0=1\subset U_1\subset U_2\subset A, $$ be the corresponding
series of ideals.   Take $H_1=G\bigcap(1+U_1)$ which is the
kernel of the action $G$ in $A/U_1$.  This group is locally
nilpotent. In $U_1$ there is a directed system of the nilpotent
ideals $U_\alpha$ of the algebra $A$. All
$G\bigcap(1+U_\alpha)=H_\alpha$ are the nilpotent normal
subgroups in $G$ and they constitute a directed system which
generates $H_1$.

Take, further,  $G\bigcap(1+U_2)=H_2.$ This is a locally
nilpotent normal subgroup in $G$ which coincides the kernel of
the action $G$ in $A/U_2$. The group $H_1$ is the kernel of
action of the group $H_2$ in $A/U_1$. This action is
unitriangular and $H_2/H_1$ is a nilpotent group.

Consider a representation $G\to A/U_2$. It corresponds the
faithful representation $G/H_2\to A/U_2$. There is also an
embedding $A/U_2\to M_n(K)$. This implies the faithful embedding
$G/H_2\to \GL_n(K)$.
$\square$

Observe also the following. Let $K=\sum_\alpha P_\alpha$. For
every $P_\alpha$ take an ideal $U_\alpha$ with
$K/U_\alpha\backsimeq P_\alpha$ and $\bigcap_\alpha=0$. For every
$\alpha$ there is a homomorphism $\GL_n(K)\to \GL_n(P_\alpha)$.
Its kernel is the congruence-subgroup in $\GL_n(K)$ modulo the
ideal $U_\alpha$. This leads to the presentation of $GL_n(K)$ as
a subdirect product of the groups $\GL_n(P_\alpha)$.


 \section{ Theorems 2 and 3}

 Let us repeat the formulations:

 1. Every nil-$PI$-group is locally nilpotent -(Theorem 2).

 2. In every $PI$-group the set of all nil-elements
 coincides with the locally  nilpotent radical $\HP(G)$ - (Theorem
 3).


  Theorem 2
   follows from Theorem 3.
  Indeed, if every element $g$ is a nil-element then $G=HP(G)$ and
  $G$ is locally nilpotent.

  As for Theorem 3
  it is, in fact, proved in
  \cite{Pla}, \cite{To1}, \cite{To2}, \cite{Pl3}. For the sake of the  self-completeness
  of the text we will give here a  proof of Theorem 3. 
  Another reason is that the same scheme works for the proof of
  Theorem 6.
   We split the proof for 3 steps.

  1. Let, first, $G$ be a linear group, i.e., $G\subset \GL_n(P)$.
  Check that every its nilpotent element lies in the radical
  $\HP(G)$.

  Let $H$ be a subgroup in $G$ generated by all its nil-elements.
  Show that $H$ is locally solvable. Take in $H$ a subgroup $H_0$
  which is generated by a finite number of nil-elements. According
  to well-known theorem of A.~Malcev \cite{Ma} there is a system of normal
  subgroups $T_\alpha$, $\alpha\in I$ in $H_0$ with the trivial
  intersection and with the finite quotients
  $H_\alpha=H_0/T_\alpha.$ These quotients $H_\alpha$ are linear groups of the same
  dimension $n$ over finite fields. Every $H_\alpha$ is generated by nil-elements
  and thus nilpotent by Baer's theorem. Therefore, $H_\alpha$ is solvable.
  Observe that all these $H_\alpha$ has the solvable length bounded
  by the number which depends only on $n$. Then $H_0$ is also
  solvable and $H$ is locally solvable. It is known that locally
  solvable linear group is solvable \cite{Su}.
Thus, $H$ is a solvable normal subgroup in $G$ generated by
nil-elements. According to \cite{Pl1} such a group is locally
nilpotent. All nil-elements of the group $G$ lie in $\HP$-radical
of $G$.

2. Consider further the case $G\subset \GL_n(K)$. The group $G$
is approximated by subgroups of linear groups $ \GL_n(P_\alpha)$.
As before, let $H$ be the subgroup in $G$ generated by all
nil-elements. This $H$ is approximated by subgroups
$H_\alpha\subset \GL_n(P_\alpha)$. The subgroups $H_\alpha$ are
generated by nil-elements and, hence, are solvable. The derived
lengths are bounded for all $H_\alpha$. Therefore, $H$ is
solvable. Since $H$ is generated by nil-elements, then $H$ is
locally nilpotent invariant subgroup. Every nil-element lies in
$H$, and therefore, in $\HP(G).$

3. General case. We have a chain $$ H_0=1\subset H_1\subset
H_2\subset H_3\subset G,$$ where $H_1$ and $H_2$ are the same as
in Proposition 3.1., and $H_3/H_2=HP(G/H_2)$. Let $g$ be a
nil-element in $G$. Take a nil-element $\bar g= gH_2$ in
$G/H_2$.  We have $\bar g\in H_3/H_2$, $g\in H_3$. Using again
\cite{Pl1}, we have $g\in R(H_3)$. Since $R(H_3)\subset R(G)$,
then $g\in R(G). \ \square $


As we have mentioned, Theorem 2 is close to Procesi-Tokarenko
theorem on periodic $PI$-groups. In some sense the theorem 2 is
related also to the theorem from \cite{WZ} where the profinite
completion of a residually finite group is considered. The
profinite setting also allows to proceed from nil-groups (not from
Engel groups).

\section{Generalizations}

Let $u=u(x,y)$ denote the elements of the free group $F_2=F(x,y)$.
Consider a sequence $\overrightarrow{u} = u_1, u_2, u_3, \ldots$.
Such a sequence is called correct, if

1. $u_n(a,1)=1$ and $u_n(1,g)=1$ for every n, every group $G$ and
every elements $a,g\in G$.

2. If $u_n(a,g)=1$ then for every $m>n$ we have $u_m(a,g)=1$ where
$a,g$ are the elements from $G.$

Thus, if the identity $u_n(x,y)\equiv 1$ is fulfilled in $G$ then
for every $m>n$ the identity  $u_m(x,y)\equiv 1$ also holds in
$G$.

For every correct sequence $\overrightarrow{u}$ consider the class
of groups $\Theta=\Theta (\overrightarrow{u})$ defined by the
rule: $G\in \Theta$ if there exists $n$ such that the identity
$u_n(x,y)\equiv 1$ holds in $G$.

For every group $G$ denote by $G(\overrightarrow{u})$ the subset
in $G$ defined by the rule: $g\in G(\overrightarrow{u})$, if for
every $a\in G$ there exists $n=n(a,g)$ such that $u_n(a,g)=1$. The
elements of $G(\overrightarrow{u})$ are viewed  as Engel elements
in respect to the given correct sequence $\overrightarrow{u}$. We
used to call these elements as $\overrightarrow{u}-$Engel-like
elements.

If $\overrightarrow{u}=\overrightarrow{e}=e_1,\ldots,e_n $ where
the words $e_n(x,y)$ are defined by
$$e_1(x,y)=[x,y],\ldots,e_n(x,y)=[e_{n-1}(x,y),y],\ldots,$$
 then
$\Theta(\overrightarrow{e})$ is the class of all Engel groups. In
case of finite groups the class $\Theta(\overrightarrow{e})$
coincides with the class of finite nilpotent groups.

For finite groups $G$ the set $G(\overrightarrow{e})$ coincides
with the nilpotent radical of the group $G$.

\begin{prob} Describe $\overrightarrow{u}$ such that
$\Theta(\overrightarrow{u})$ is the class of finite solvable
groups.
\end{prob}

Concerning this problem see \cite{BGGKPP}, \cite{Wi}.

\begin{prob} Construct a sequence $\overrightarrow{u}$ such that
$G(\overrightarrow{u})$ is the solvable radical of every finite
group $G$.
\end{prob}

It is not known whether there exist such $\overrightarrow{u}$. In
both problems above emphasis is made on the fact that we are
looking for two variable sequences. However, the similar problems
can be considered also in the general case were the number of
variables is not restricted.

In particular, we will consider some other approach to the problem
of solvable radical description, which makes sense for finite
groups too.

\section{Further generalizations}

For each given correct sequence $\overrightarrow{u}$ define a new
set  $\overrightarrow{\overrightarrow{u}}$. Consider the free
group $F=F(X,y)$ where $X=\{x_1,x_2,\ldots,x_k,\ldots\}$, and $y$
is a distinguished variable. We will index words from $F$ by the
sequences of natural numbers $(n_1,n_2,\ldots,n_k)$. Define the
words
$$
u_{(n_1,n_2,\ldots,n_k)}(x_1,x_2,\ldots,x_k;y)
$$
 by the rule: $u_{n_1}(x_1,y)$ coincides with the corresponding
 element of the sequence $\overrightarrow{u}$. Then, by induction,
$$
u_{(n_1,n_2,\ldots,n_k)}(x_1,x_2,\ldots,x_k;y)=
u_{n_k}(x_k;u_{(n_1,n_2,\ldots,n_{k-1})}(x_1,x_2,\ldots,x_{k-1};y)).
$$
The considered words obtained by superposition of two variable
words.

It is easy to see that the following associativity takes place:
$$ u_{(n_1,n_2,\ldots,n_k)}(x_1,x_2,\ldots,x_k;y)= $$ $$
u_{(n_{l+1},. \ldots,n_k)}(x_{l+1},\ldots,x_k;
u_{(n_1,n_2,\ldots,n_l)}(x_1,x_2,\ldots,x_l;y)). $$ In
particular, $$ u_{(n_1,n_2,\ldots,n_k)}(x_1,x_2,\ldots,x_k;y)= $$
$$ u_{(n_2,\ldots,n_k)}(x_2,\ldots,x_k;u_{n_1}(x_1;y)) $$
Correctness of the initial sequence $\overrightarrow{u}$ induces
some correctness of the system
$\overrightarrow{\overrightarrow{u}}.$

For example, if for $l<k$ the group $G$ satisfies the identity $$
u_{(n_l,\ldots,n_k)}(x_l,\ldots,x_k;y)\equiv 1, $$ or the identity
$$
 u_{(n_1,\ldots,n_{l-1})}(x_1,\ldots,x_{l-1};y)\equiv 1,
  $$
then $G$ satisfies the identity $$
u_{(n_1,\ldots,n_k)}(x_1,\ldots,x_k;y)\equiv 1. $$ There are also
other relations of such kind.

For the given system $\overrightarrow{\overrightarrow{u}}$
consider the class of groups
$\Theta=\Theta(\overrightarrow{\overrightarrow{u}}).$ By
definition, a group $G$ belongs to $\Theta$ if an identity of the
form
$$
u_{(n_1,\ldots,n_k)}(x_1,\ldots,x_k;y)\equiv 1.
$$
holds in $G$. From the observations above follows that the class
$\Theta$ is a pseudovariety of groups. Besides that, for every
group $G$  we define a class of elements
$G(\overrightarrow{\overrightarrow{u}})$ by the rule: $g\in
G(\overrightarrow{\overrightarrow{u}})$ if for some $k=k(g)$ and
for every sequence $(a_1,\ldots,a_k)$ of elements in $G$ there is
a set $(n_1,\ldots, n_k)$ such that
$$
u_{(n_1,\ldots,n_k)}(a_1,\ldots,a_k;g)= 1.
$$
is fulfilled. Here, the set $(n_1,\ldots, n_k)$ should be
compatible with the set $(a_1,\ldots,a_k). $ This means that $n_1$
depends on $a_1$ and $g$, and does not depend on
$(a_2,\ldots,a_k)$; $n_2$ depends on $a_1,a_2,$ and $g$, and does
not depend on $(a_3,\ldots,a_k)$, etc.; $n_k$ depends on
$a_1,a_2,\ldots,a_k$ and $g$.

Here arises a general problem of some description of the sets
$G(\overrightarrow{\overrightarrow{u}})$ for the different
$\overrightarrow{u}$.

Consider a special case of
$G(\overrightarrow{\overrightarrow{u}})$. Denote  $\varepsilon=
\overrightarrow{e}$ and take the sequence
$$
\varepsilon:e_1,e_2,\ldots, e_n, \ldots
$$
Consider the system $\overrightarrow{\varepsilon}$ and using this
system define {\it quasi-nil elements} in groups. An element $g\in
G$ is called {\it quasi-nil} if $g\in
G(\overrightarrow{\varepsilon}).$ This means that for $g$ there is
$k=k(g)$, such that for any sequence $a_1,\ldots, a_k$, $a_i\in
G$,  there is a compatible set $(n_1,\ldots,n_k)$ such that
$$
\varepsilon_{(n_1,\ldots,n_k)}(a_1,\ldots,a_k;g)= 1.
$$

For the sequence $\varepsilon$ we have also the class of groups
$\Theta(\overrightarrow{\varepsilon})$. The groups from this class
can be considered simultaneously as generalized nilpotent and
generalized solvable groups.

Denote by $E_{(n_1,\ldots,n_k)}$ the variety defined by the
identity
$$
\varepsilon_{(n_1,\ldots,n_k)}(x_1,\ldots,x_k;y)\equiv 1.
$$
The class $\Theta$ is the union of such varieties. The variety of
the type $E_{(1,\ldots,1)}$ is the nilpotent variety, while the
variety of the type $E_{(1,2,\ldots,2)}$ contains the solvable
subvariety. Besides that a product of varieties of the type
$E_{(n_1,\ldots,n_k)}$  is a subvariety in the variety of the same
type. This observation relates, in particular, to the product
$E_{n_1}E_{n_2}\cdots E_{n_k}$. This variety lies in the variety
$E_{(n_k, n_{k-1}+1,\ldots,n_1+1)}$.

Return now to quasi-nil elements in groups. Let $k=k(g)$ be the
minimal number such that for every $(a_1,\ldots,a_k)$ there is a
compatible set $(n_1,\ldots,n_k)$ with
$$
\varepsilon_{(n_1,\ldots,n_k)}(a_1,\ldots,a_k;g)=1.
$$

We call such $k=k(g)$ the nil-order of $g$. Nil-order 1 means that
the element is a nil-element, nil-order 2 means that the element
is not nil, but for $a_1$ and $a_2$ there are $n_1$, $n_2$ with
$\varepsilon _{(n_1,n_2)}(a_1,a_2;g)=1$. In general  for $k-1$ we
have some elements $(a_1^0,\ldots,a_{k-1}^0)$ such that

$$
\varepsilon_{(n_1,\ldots,n_{k-1})}(a^0_1,\ldots,a^0_{k-1};g)\neq
1.
$$
for the arbitrary compatible $(n_1,\ldots,n_{k-1})$.

Let us add to $(a_1^0,\ldots,a_{k-1}^0)$ an arbitrary element $a$.
Then for the sequence $(a_1^0,\ldots,a_{k-1}^0,a)$ there is a
corresponding set $(n_1^0,\ldots,n_{k-1}^0,n)$ with the condition
$$ e_n(a;g_0)=1, $$
 where $$
g_0=\varepsilon_{(n_1^0,\ldots,n_{k-1}^0)}(a_1^0,\ldots,a_{k-1}^0;g).
$$ Here, the element $g_0$ is not trivial, the element $a$ does
not depend on $g_0$. The equality $e_n(a,g_0)=1$ now means that
the element $g_0 $ is a non-trivial nil-element.

Simultaneously, we proved the following

\begin{prop} If a group $G$ contains a non-trivial quasi-nil element $g$
then $G$ contains also a non-trivial  nil-element $g_0$.
\end{prop}

Note now the next two properties related to the definition of the
quasi-nil element.

1. Let $H$ be a subgroup in $G$ and $g\in H$ be a quasi-nil in
$G$. Then $g$ is quasi-nil element in $H$.

2. Let a surjection $\mu: G\to H$ be given and let $g$ be a
quasi-nil element in $G$. Then $\mu(g)$ is a quasi-nil element in
$H$.

Indeed, take $k=k(g)$ and the corresponding presentation
$$
\varepsilon_{(n_1,\ldots,n_{k})}(a_1,\ldots,a_{k};g)= 1.
$$
Then
$$
\varepsilon_{(n_1,\ldots,n_{k})}(\mu(a_1),\ldots,\mu(a_{k};\mu(g))=
1.
$$

Here, $\mu(a_1), \ldots, \mu(a_n)$ are arbitrary elements in $H$.

It is clear that along with quasi-nil elements it is quite natural
to define quasi-Engel elements which generalize Engel elements.

\section{Some radicals}

Let $G$ be a group. Consider in $G$ the locally nilpotent radical
$\HP(G)=R(G)$ and the locally noetherian radical ${NR}(G)$. The
corresponding upper radicals will be denoted by  $\widetilde
{R}(G)$ and $\widetilde{NR}(G)$. These radicals are obtained by
iterations of the initial $R(G)$ and  ${NR}(G).$ Namely, consider
the series (upper radical series) $$ 1=R_0\subset
R=R_1\subset\ldots\subset R_\alpha\subset\ldots , $$ where
$R_{\alpha+1}/R_\alpha$ is $R(G/R_\alpha)$. Such a series
terminates at some $R_\gamma=\widetilde R(G)$. Then $\widetilde
R(G)$ is the upper radical for the radical $R(G)$. The factor
group $G/\widetilde R(G)$ is locally nilpotent  semi-simple. i.e.,
it does not contain non-trivial locally nilpotent normal
subgroups.

The radical $\widetilde R$ is defined also by the class of groups
$G$ which has ascending normal series with locally nilpotent
factors. Such groups are called radical groups (see \cite{Pl1}).
In finite groups the radical $\widetilde R(G)$  coincides with the
solvable radical of a group.

The radical $\widetilde{NR}(G)$ is defined following the same
scheme as  for the radical $\widetilde R(G)$. If
$\widetilde{NR}(G)=G$, the group $G$ is called noetherian radical
group.

\section{Theorems on radical characterization}

Let us take in the upper radical series of a group $G$ the members
with finite indexes
$$
1=R_0\subset R_1\subset\ldots R_k\subset\ldots ,
$$

\begin{prop} An element $g$ which belongs $R_k$ for some $k$ and
does not belong to $R_{k-1}$ is a quasi-nil element of the
nil-order $k$.
\end{prop}

Proof. For the case $g \in R_1$ this is true. Further we use
induction. Suppose that for $g \in R_{k-1}$ it is proved that the
nil-order of this $g$ is $\leqslant k-1$. Let $g \in R_k$. Take a
sequence of elements $a_1, \ldots, a_k$ in $G$ and for $a_1$ and
$g$ find $n_1$ with $e_{n_1} (a_1, g) \in R_{k-1}$. Apply
induction to the element $e_{n_1} (a, g)$. We have:

$$
\varepsilon_{(n_2, \ldots, n_k)} (a_2, \ldots, a_k; e_{n_1} (a_1,
g)) = 1 = \varepsilon_{(n_1, \ldots, n_k)} (a_1, \ldots, a_k; g).
$$

Hence, the nil-order of the element $g$ is $\leqslant k$. Prove
further that it is exactly $k$. Let $g$ is of the order
$l\leqslant k$. Take $a_1, \ldots, a_l, n_1, \ldots, n_l$ such
that
$$
\varepsilon_{(n_1, \ldots, n_l)} (a_1, \ldots, a_l; g) = 1 =
e_{n_l} (a_l; \varepsilon_{(n_1, \ldots, n_{l-1})} (a_1, \ldots,
a_{l-1}; g).
$$

The element $a_l$ does not depend on  $g_0=\varepsilon_{(n_1,
\ldots, n_{l-1})} (a_1, \ldots, a_{l-1}; g),$ and all these $g_0$
are nil-elements (for all $a_1,\ldots,a_{l-1}$). Some of $g_0$ are
non-trivial and all of them lie in $R_1$. Consider $G/R_1$. Here
all $\bar g_0$ are trivial and the nil-order of $\bar g$ is
$\leqslant l-1$. By the assumption of induction $\bar g\in
R_l/R_1$. Then $g\in R_l$. By the condition $g$ does not belong to
$R_{k-1}$. Then $l=k$. $\square$

\begin{prop} Let $\widetilde {NR}(G)=G$ and $\widetilde R(G)$ be the
radical. Then every quasi-nil element $g\in G$ belongs to
$\widetilde R(G)$.
\end{prop}

Proof. Let $g$ be a quasi-nil element which does not belong to
$\widetilde R(G)$. Element $\bar g= g\widetilde R(G)$ is quasi-nil
in the semi-simple group $\bar G= G/\widetilde R(G).$ If $g\neq 1$
then there exists a non-trivial nil element in $\bar G$. We came o
contradiction with the semi-simplicity of $G$.$\square$

\begin{thm} Let $\widetilde {NR}(G)=G$ and let the upper radical
series in $G$ has a finite length. Then $\widetilde R(G)$
coincides with the set of all quasi-nil elements in $G$.
\end{thm}

Proof. The proof follows from Proposition 1 and Proposition 2.

 We have
 seen that every quasi-nil element in $\widetilde {NR}(G)=G$  lies
 in $\widetilde R(G)$ for  upper radical series of any length.
 However, in this general situation we cannot state that every
 element from $\widetilde R(G)$ is quasi-nil. In order to
 include this case in the general setting we define unbounded
 quasi-nil elements. In this unbounded approach we do not fix
 $k=k(g)$, since we do not know in advance what are the length of
 words which are related with the given $g$. Thus we consider
 infinite sequences $\bar a=(a_1,a_2,\ldots,a_k,\ldots)$. We call
 an element unbounded quasi-nil, if for any $\bar a$ there are
 $k=k(\bar a, g)$ and compatible $(n_1,n_2,\ldots, n_k)$ such that
 $$
 \varepsilon_{(n_1,n_2,\ldots,n_k)}(a_1,a_2,\ldots,a_k;g)=1.
 $$

\begin{thm} For any group $G$ every element in $\widetilde R(G)$
is an unbounded quasi-nil element. \end{thm}

Proof. We start from the upper radical series

$$ 1=R_0\subset R_1\subset R_2\subset\ldots\subset
R_\alpha\subset\ldots\subset R_\gamma=R $$
and use the induction.
For $g\in R_1$ the statement is evident and let for all $\beta <
\alpha $ the statement is true. Show that every $g\in R_\alpha $
is unbounded quasi-nil.

If $\alpha$ is terminal then $g\in R_\beta$ with $\beta < \alpha$
and $g$ is unbounded quasi-nil.

Suppose now that there exists $\alpha-1$. For the given $g\in
R_\alpha$ take a sequence
$$
\bar a= (a_1, a_2, \ldots, a_k, \ldots   ).
$$
For $a_1$ we find $n_1$ with $e_{n_1}(a_1,g)\in R_{\alpha-1}$. The
element $e_{n_1}(a_1,g)$ is unbounded quasi-nil. For this element
take the sequence $a_2, \ldots, a_k, \ldots   $,  and let $(n_2,
\ldots, n_k,)$ be defined for this sequence. Here, $k$ depends
also on $a_1$.  We have
$$
\varepsilon_{(n_2,\ldots,n_{k})}(a_2,\ldots,a_{k};e_{n_1}(a_1,g))=1
$$
$$
=\varepsilon_{(n_1,n_2\ldots,n_{k})}(a_1, a_2,\ldots,a_{k};g).
$$
Thus, the element $g$ satisfies the condition to be unbounded
quasi-nil.

\section{Again about PI-groups}

Theorem 3 can be applied to finite groups. It can be also applied
to linear groups over fields and, as we will see soon, to any
$PI$-groups. In these cases the conditions of the type $\widetilde
{NR}(G)=G$ are not necessary.

\begin{thm} If $G$ is a $PI$-group, then its "solvable"
radical $\widetilde R(G)$ coincides with the set of all quasi-nil
elements.
\end{thm}

The proof of this theorem follows the scheme used in the proof of
Theorem 2. The only observation has to be taken into account is
the fact that in every solvable group its solvable radical
coincides with the set of its nil-elements. $\square$

In particular, we can state that a $PI$-group is "solvable" (in
the sense that $\widetilde R(G)=G$) if and only if all elements in
$G$ are quasi-nil.

From the theorem 2 follows that if in $PI$-group $G$ every two
elements generate a nilpotent subgroup then the whole group is
locally nilpotent. Now we consider the case when every two
elements generate a solvable subgroup.

\begin{thm} Let $G$ be a $PI$-group and let every two elements in
$G$ generate a solvable subgroup. Then $G$ is solvable modulo
locally nilpotent radical $\HP(G)$. \end{thm}

Proof. It is sufficient to consider the case when $G$ is a
subgroup in a $\GL_n(K)$, where $K$ is a direct sum of fields. If
$K$ is a field the proof follows from \cite{Th}, \cite{Su}. The
proof for the general case imitates the reduction to the field
case in the previous theorem. $\square$

\begin{remark} In every $PI$-group $G$ the group $\widetilde
{R}(G)/R(G)$ is solvable. However, it is not clear whether the
group $\widetilde {R}(G)$ is always locally solvable. \end{remark}

 \section{Conclusion}

All above can be applied to finite groups. However, the problems 1
and 2 remains open. Their solution should use the "subtle" theory
of finite groups. This is the classification of finite simple
groups and their automorphisms, equations in finite simple groups,
etc. Here some algebraic geometry can be used. Besides that, along
with Engel-like elements the corresponding Engel-like
automorphisms should be considered.


\end{document}

%% file: cd.tex
\catcode`\@=11

 \def\AMSTeXfeatures{\Plainheads 
   \let\current@vert=\AMS@vert}

 \def\Plainheads{\sh@ftdiam=0.05em
   \getlabeldims
   \let\vshaftfill=\plnvsolidfill
   \let\hshaftfill=\plnhsolidfill
   \let\th@rhead=\plnrhead
   \let\th@lhead=\plnlhead
   \let\th@dnhead=\plndnhead
   \let\th@uphead=\plnuphead}
 
 \def\glet{\global\let}

 \def\LaTeXfeatures{\catcode`\@=11
   \ifx\@clnwd\undefined \nol@g
      \input ltxcode.tex \dol@g \fi
   \ltxheads \let\current@vert=\new@vert
   \providelto \catcode`\@=\active}

 \def\nol@g{\def\wlog{\edef\garbage}}
 \def\dol@g{\let\wlog=\wl@g} \let\wl@g=\wlog
 \nol@g 

 \newbox\ltobox
 \def\providelto{{\setbox\z@=
   \hbox{$\to$}\minharrlen=\wd\z@
   \global\setbox\ltobox=\hbox{$\activeat>>>$}}
   \def\lto{\mathrel{\copy\ltobox}}}

 \def\ltxheads{\sh@ftdiam=\@wholewidth
   \getlabeldims
   \let\vshaftfill= \ltxvsolidfill
   \let\hshaftfill=\ltxhsolidfill
   \let\th@rhead=\ltxrhead
   \let\th@lhead=\ltxlhead
   \let\th@dnhead=\ltxdnhead
   \let\th@uphead=\ltxuphead}
 {\catcode`\@=\active
   \gdef@#1{\csname #1\string@at\endcsname}
   \glet\activeat=@}
 \def\def@#1{\expandafter\def\csname #1@at\endcsname}

 \def@>#1>#2>{\@rrow R{#1}{#2}}
 \def@<#1<#2<{\@rrow L{#1}{#2}}
 \def@ V#1V#2V{\@rrow V{#1}{#2}}
 \def@ A#1A#2A{\@rrow A{#1}{#2}}
 \def@/#1/#2/#3/{\@rrow{#1}{#2}{#3}}
 \def@.{\ifodd\row\ifmmode\noharrow
     \else\leavevmode.\spacefactor3000 \fi
   \else\novarrow\fi}
 \def@={\ifodd\row\harrow\hequalfill{}{}%
   \else\varrow\vequalfill{}{}\fi}
 \def@:#1{\ifx=#1\harrow\deffill{}{}%
   \else\leavevmode\null:#1\fi}
 \def@|{\current@vert}
  \def\AMS@vert{\varrow\vequalfill{}{}}
  \def\new@vert#1|#2|{\ifodd\row
   \let\nextarrow\vertexvarrow
   \else\let\nextarrow\varrow\fi
   \nextarrow\vshaftfill{#1}{#2}}
 \def@-{\ifmmode\let\next\hl@ne
   \else\let\next\AMSatdash \fi \next}
  \def\hl@ne#1-#2-{\harrow\hshaftfill{#1}{#2}}
  \def\AMSatdash{\let\next\relax\leavevmode
    \def\next@{\ifx\next-%
      \def\next-{\futurelet\next\nextii@}%
     \else\def\next{\hbox{-}}\fi\next}%
    \def\nextii@{\ifx\next-\def\next-{\hbox{---}}%
      \else\def\next{\hbox{--}}\fi\next}%
    \futurelet\next\next@}
 \def@(#1){\tweenarrows{#1}}
 \def@[#1]{\setsp@n#1\relax\activeat}
 \def\fiberbox{\hbox{$\vcenter{\hr@le\hbox{\vr@le
   \kern1ex\vbox{\kern1.2ex}\vr@le}\hr@le}$}}
  \def\hr@le{\hrule height \sh@ftdiam}
  \def\vr@le{\vrule width \sh@ftdiam}
 \def@+#1+#2+#3+{\ifodd\row \harrow{#1}{#2}{#3}%
   \else \varrow{#1}{#2}{#3}\fi}

 \def\Dnarrfill{\vequalfill\Dnhe@d}
 \def\Uparrfill{\Uphe@d\vequalfill}
 
 \def\ontofill{\rtarrfill\kern-0.3em 
   \th@rhead\kern 0.3em} 

 \def\rtarrfill{\hshaftfill\th@rhead}
 \def\ltarrfill{\th@lhead\hshaftfill}
 \def\dnarrfill{\vshaftfill\th@dnhead}
 \def\uparrfill{\th@uphead\vshaftfill}
 \def\hequalfill{\plnhfill=}
 \def\deffill{:\plnhfill=}
 \def\plnvextfill#1{\setbox\z@
   \hbox{\the\textfont3 #1}%
   \dimen@=\dp\z@\advance\dimen@\ht\z@
   \copy\z@ \kern-\dimen@ 
   \cleaders\copy\z@ \vfill
   \kern-\dimen@ 
   \box\z@}
 \def\plnhfill#1{$\m@th\mkern-1.5mu\mathord#1\mkern-6mu
    \cleaders\hbox{$\mkern-2mu\mathord#1\mkern-2mu$}\hfill
    \mkern-6mu\mathord#1\mkern-1.5mu$}
 \def\vequalfill{\plnvextfill{\char'167}}
 \def\plnvsolidfill{\plnvextfill{\char'077}}
 \def\plnhsolidfill{\plnhfill-}
 \def\ltxhsolidfill{\leaders\hrule height\topofshaft depth\botofshaft
   \hfill}
 \def\ltxvsolidfill{\leaders\vrule width\sh@ftdiam\vfill}
 \def\hdashfill{\hd@sh\wd@sh
   \xleaders \hbox{\wd@sh\hd@sh\wd@sh}\hfill
   \wd@sh\hd@sh}
 \def\vdashfill{\vd@sh\wd@sh
   \xleaders \vbox{\wd@sh\vd@sh\wd@sh}\vfill
   \wd@sh\vd@sh}
 \def\dashed{\ifinmeasureCD\else
    \ifodd\row\option{\let\hshaftfill=\hdashfill}%
   \else\option{\let\vshaftfill=\vdashfill}\fi\fi}

 \newdimen\CDstrutht  \newdimen\CDstrutdp
 \newdimen\CDstrutlen \CDstrutlen=\CDstrutht
   \advance\CDstrutlen by \CDstrutdp

 \def\CDstrut{\vrule
   height \ifnum\row=1 \z@\else\CDstrutht \fi
   depth \ifnum\row=\numrows \z@ \else\CDstrutdp \fi
   width\z@}

 \newdimen\CDarrsurr \CDarrsurr=0.375em
 \newdimen\CDdashlen
 \newdimen\CDvarrlen \CDvarrlen=1.5\baselineskip
 \newdimen\minharrlen 
  \setbox\z@\hbox{$\longrightarrow$} \minharrlen=\wd\z@
 \newdimen\minCDharrlen \minCDharrlen=2.5em 
\newdimen \minc@lwd
\def\findminc@lwd{\minc@lwd=2\CDarrsurr
  \advance\minc@lwd\minCDharrlen}

 \newdimen\sh@ftdiam

 \newdimen\labelsurr \labelsurr=1.25 em

\newcount\sp@ncnt \sp@ncnt=\@ne
\newcount\sp@ncnt@ \sp@ncnt@=\@ne
\newdimen\@rrwd \newdimen\@rrdp

 \def\adjustbot#1{\option{\advance\@rrdp#1\relax}}
 
\def\pushvertex#1{\global\p@shlen#1\relax
   \global\let\maybepush=\dopush}

 \newdimen\p@shlen \p@shlen=\z@

 \let\maybepush=\relax
 \def\dopush{\ifinmeasureCD 
   \advance\locdimen by -\p@shlen 
   \else\advance \@rrwd by -\p@shlen \fi 
   \global\let\maybepush=\relax \global\p@shlen=\z@\relax}

 \def\span@ne{\global\sp@ncnt=\@ne\relax}
 \def\setsp@n#1#2{\global\sp@ncnt=#1\relax
   \ifx\relax#2\relax\else\global\sp@ncnt@=#2\relax\fi}

 \def\plnrhead{\llap{$\rightarrow\mkern-1.5mu$}}
 \def\plnlhead{\rlap{$\mkern-1.5mu\leftarrow$}}

 \def\clap#1{\hbox to \z@{\hss #1\hss}}

 \def\plndnhead{\hbox{\the\textfont3 \char'171}}
 \def\plnuphead{\hbox{\the\textfont3 \char'170}}
 \def\Dnhe@d{\hbox{\the\textfont3 \char'177}}
 \def\Uphe@d{\hbox{\the\textfont3 \char'176}}

 \def\ltxrhead{\raise\@xisheight
   \llap{\smash{\@linefnt\@getrarrow(1,0)}}}
 \def\ltxlhead{\raise\@xisheight
   \rlap{\@linefnt\@getlarrow(-1,0)}}
 \def\ltxuphead{\setbox\z@=\rlap{%
   \kern\@halfwidth\@linefnt\char'66}%
   \copy\z@\kern-\ht\z@}
 \def\ltxdnhead{\setbox\z@=\rlap{%
   \kern\@halfwidth\@linefnt\char'77}%
   \ht\z@=\z@\box\z@}

 \def\wd@sh{\kern0.5\CDdashlen}
 \def\hd@sh{\vrule height\topofshaft depth\botofshaft
    width\CDdashlen}
 \def\vd@sh{\hrule height\CDdashlen
   depth\z@ width\sh@ftdiam}

\def\xylist{14{3434}13{2414}12{1723}%
  23{1413}34{1153}11{0867}43{0707}%
  32{0580}21{0414}31{0291}41{0}}
\newcount\tgtcnt@
\def\find@xyargs{\dimen@=\@rrdp
  \advance\dimen@ by \CDstrutlen
  \tgtcnt@=\dimen@ \dimen@=\@rrwd 
  \divide\dimen@ by \@m 
  \divide \tgtcnt@ by \dimen@ 
  \expandafter\testxy\xylist\relax
  \unitlength=\@xarg\@rrdp
  \divide\unitlength by\@yarg\relax}
\def\testxy#1#2#3{\ifnum\tgtcnt@>#3
    \@xarg=#1\relax \@yarg=#2\relax
    \let\next=\ignorerest
  \else\let\next\testxy\fi\next}
\def\ignorerest#1\relax{\relax}

\let\scalefactor=\@ne
\def\SWarrow{\find@xyargs\vector
  (-\@xarg,-\@yarg)\scalefactor\hskip-\wd\@linechar}
\def\NWarrow{\find@xyargs\vector
  (-\@xarg,\@yarg)\scalefactor\hskip-\wd\@linechar}
\def\NEarrow{\find@xyargs\vector
  (\@xarg,\@yarg)\scalefactor}
\def\SEarrow{\find@xyargs\vector
  (\@xarg,-\@yarg)\scalefactor}
\def\rightupline{\find@xyargs\@linelen=\scalefactor
     \unitlength\@sline}
\def\rightdownline{\find@xyargs\@yarg=-\@yarg\relax
     \@linelen=\scalefactor\unitlength\@sline}

\def\Sim{\ifodd\row\setbox\z@=\hbox{$\sim$}\dimen@=\ht\z@
 \advance\dimen@ by -\@xisheight
  \vbox{\box\z@\kern-\@xisheight\kern\dimen@}%
  \else\hbox{$\wr$}\fi}

\def\harrow#1#2#3{\inmeasureCDtrue\findminarrwd
  {#2}{#3}{\sp@ncnt\minharrlen}\inmeasureCDfalse\span@ne
  \mathrel{\hbox{\options\hplace{#1}\ulabel{#2}\dlabel{#3}}}}

\def\noharrow{\harrow\hfill{}{}}
\def\vertexvarrow#1#2#3{\findarrdp \@rrwd=\z@ \setsp@n\@ne\@ne
  \vbox to \z@{\kern-1.2\CDstrutht
  \rlap{\options\vplace{#1}\llabel{#2}\rlabel{#3}}\vss}}

\newif\ifinmeasureCD
\def\measurelabel#1{\setbox\z@
  \hbox{$\scriptstyle#1\kern\labelsurr$}%
  \ifdim\wd\z@>\@rrwd \@rrwd=\wd\z@\fi}
\def\findminarrwd#1#2#3{\@rrwd=#3\relax
   \measurelabel{#1}\measurelabel{#2}}
\def\findCDarrwd#1#2{\@rrwd=\minCDharrlen
   \measurelabel{#1}\measurelabel{#2}%
  }

\newcount\row \row=\@ne \newcount\col \col=\@ne 
 \newcount\numrows 
\numrows=\@ne
 \newcount\numcols
\newcount\arrspan \newdimen\vrtxhalfwd  \newbox\tempbox

\def\DANABUG{\advance\col by \@ne
 \@rrwd=\minCDharrlen
  \advance\@rrwd by \vrtxhalfwd
  \advance\@rrwd by \CDarrsurr
  \ifnum\col>\numcols \numcols=\col
     \newlocdimen{col\the\col}\locdimen=\@rrwd 
  \else \ifdim\@rrwd>\c@l \c@l=\@rrwd\fi\fi}

\def\drop#1\\{
  \findvrtxhalfsum\DANABUG\advance\row by 2 \measureinit}

\def\measureinit{\col=\@ne \vrtxhalfwd=-\CDarrsurr\arrspan=\@ne\@rrwd=\z@
   \setbox\tempbox=\hbox\bgroup$}
\def\measure{
  \let\harrow\measureCDarrow
  \let\CDCR=\measureCR 
   \findminc@lwd 
  \inmeasureCDtrue
  \row=\@ne \numcols=\z@ \measureinit}

\def\endmeasure{\findvrtxhalfsum\DANABUG
  \numrows=\row 
  \inmeasureCDfalse}

\def\newlocdimen#1{\advance\dimenc@unt by \@ne
  \ifnum\dimenc@unt<\insc@unt
     \else\errmessage{No room for the CD}\fi
  \dimendef\locdimen=\dimenc@unt
  \expandafter\dimendef\csname#1\endcsname=\dimenc@unt}

 \def\r@wc@l{\csname row\the\row col\the\col\endcsname}
 \def\c@l{\csname col\the\col\endcsname}

 \def\findvrtxhalfsum{$\egroup
  \newlocdimen{row\the\row col\the\col}
  \locdimen=\vrtxhalfwd 
  \vrtxhalfwd=0.5\wd\tempbox 
  \advance\vrtxhalfwd by \CDarrsurr
  \advance\locdimen by \vrtxhalfwd 
  \advance\@rrwd by \locdimen 
  \maybepush
  \divide\@rrwd by \arrspan\relax
  \ifdim\@rrwd<\minc@lwd
    \ifnum\col>\@ne \@rrwd=\minc@lwd\fi \fi
  \loop 
    \ifnum\col>\numcols \numcols=\col
       \newlocdimen{col\the\col}
       \locdimen=\@rrwd 
    \else \ifdim\@rrwd>\c@l \c@l=\@rrwd\fi \fi
   \ifnum\arrspan>\@ne
      \advance\arrspan by -1 \advance\col by \@ne
  \repeat }

 \def\measureCDarrow#1#2#3{\findvrtxhalfsum
   \arrspan=\sp@ncnt\relax\global\sp@ncnt=1\relax
   \advance\col by \@ne
   \findCDarrwd{#2}{#3}%
   \setbox\tempbox=\hbox\bgroup$}

 \newcount\dr@tn \dr@tn=\z@
 \def\locate#1:#2{\ifinmeasureCD\else
   \count@=-#1
   \multiply\count@ by 2
   \advance\count@ by #2
   \dimen@=\count@\@rrwd
   \ifnum\dr@tn=\@ne\relax \else\dimen@=-\dimen@ \fi
   \dimen@i=\@rrdp
   \ifnum\dr@tn>\z@\advance\dimen@i by \CDstrutlen \fi
   \dimen@i=\count@\dimen@i
   \count@=#2 \multiply\count@ by 2
   \divide\dimen@ by \count@
   \divide\dimen@i by \count@
   \lift\dimen@i\nudge\dimen@\fi}

\def\betweenCDrows{\advance\row by \@ne \col=\@ne
\options}

\def\hbegin{\hbox\bgroup\kern\c@l \kern-\r@wc@l$}
\def\hend{$\glet\maybepush\relax \CDstrut\egroup}
\def\vbegin{\setbox\tempbox=\hbox\bgroup$}
\def\vend{$\egroup\ht\tempbox=\z@\dp\tempbox\CDvarrlen
  \box\tempbox}
\def\setCD{\let\harrow=\setCDarrow
  \let\CDCR=\setCR 
  \row=\@ne \col=\@ne \hbegin}
\let\endsetCD=\hend 

\def\findarrwd{\@rrwd=\z@ \count@=\col \advance\count@ by\sp@ncnt
  \loop\ifnum\count@>\col \advance\count@ by -1
      \advance\@rrwd by\csname col\the\count@\endcsname\repeat}
\def\setCDarrow#1#2#3{\kern\CDarrsurr\advance\col by \@ne
  \findarrwd \advance\@rrwd by -\r@wc@l  
  \@rrdp=\z@ 
  \maybepush
  \advance\col by -\@ne \advance\col by \sp@ncnt \span@ne
  \hbox to \@rrwd{\options
   \@rrwd=\scalefactor\@rrwd\hss
   \hplace{#1}\ulabel{#2}\dlabel{#3}\hss}%
   \kern\CDarrsurr}

\newdimen\labspacei 
\newdimen\labspaceii 

\newdimen\@xisheight
  \@xisheight=\the\fontdimen22\textfont2
\newdimen\labelskip
  \labelskip=\the\fontdimen10\textfont3 
\newdimen\topofshaft
\newdimen\botofshaft
\newdimen\botofulabel
\newdimen\topofdlabel
\def\getlabeldims{
  \topofshaft=0.5\sh@ftdiam
  \botofshaft=\topofshaft
  \advance\topofshaft by \@xisheight  
  \advance\botofshaft by -\@xisheight  
  \botofulabel=\topofshaft
  \advance\botofulabel by \labelskip
  \topofdlabel=\botofshaft
  \advance\topofdlabel by \labelskip}

\def\ulabel{\ifnum\row=\@ne\let\next\ulabeli
   \else\let\next\ulabellap\fi\next}
\def\ulabeli#1{\vbox{
  \clap{\kern-\@rrwd$\scriptstyle#1$}%
  \kern\botofulabel}\maybeoffset}
\def\ulabellap#1{\vbox to \z@{\vss
  \clap{\kern-\@rrwd$\scriptstyle#1$}%
  \kern\botofulabel}\maybeoffset}
\def\dlabel{\ifnum\row=\numrows\let\next\dlabeli
   \else\let\next\dlabellap\fi\next}
\def\dlabeli#1{\vtop{\kern\topofdlabel
  \clap{\kern-\@rrwd$\scriptstyle#1$}%
  }\maybeoffset}
\def\dlabellap#1{\vbox to \z@{\kern\topofdlabel
  \clap{\kern-\@rrwd$\scriptstyle#1$}%
  \vss}\maybeoffset}
\def\rlabel#1{\vbox to \z@{\vss
  \rlap{\kern\labelskip$\scriptstyle#1$}%
  \vss\kern-\@rrdp}\maybeoffset}
\def\llabel#1{\vbox to \z@{\vss
  \llap{$\scriptstyle#1$\kern\labelskip}%
  \vss\kern-\@rrdp}\maybeoffset}
\def\swlabel#1{\vtop{\kern0.5\@rrdp
  \llap{$\scriptstyle#1$\kern\labelskip\kern-0.5\@rrwd}
  }\maybeoffset}
\def\nwlabel#1{\vbox{
  \llap{$\scriptstyle#1$\kern\labelskip\kern-0.5\@rrwd}%
  \kern-0.5\@rrdp}\maybeoffset}
\def\selabel#1{\vtop{\kern0.5\@rrdp
  \rlap{\kern0.5\@rrwd\kern\labelskip$\scriptstyle#1$}%
  }\maybeoffset}
\def\nelabel#1{\vbox{
  \rlap{\kern0.5\@rrwd\kern\labelskip$\scriptstyle#1$}%
  \kern-0.5\@rrdp}\maybeoffset}
\def\cplace#1{\vbox to \z@{\vss
  \clap{$#1$\kern-\@rrwd}%
  \kern-\@rrdp\vss}\maybeoffset}
\def\hplace#1{\hbox to \@rrwd{#1}\maybeoffset}
\def\vplace#1{\clap{\vbox to \z@{#1\kern-\@rrdp}}\maybeoffset}

\newdimen\nudgeamount \nudgeamount=\z@
\newdimen\liftamount \liftamount=\z@
\let\maybeoffset\relax
\newbox\offsetbox \newdimen\lastheight
\def\dooffset{
  \setbox\offsetbox=\lastbox \lastheight=\ht\offsetbox 
  \setbox\offsetbox=\vbox{\kern-\liftamount\box\offsetbox}%
  \ht\offsetbox=\lastheight
  \kern\nudgeamount\box\offsetbox\kern-\nudgeamount
  \global\nudgeamount=\z@ \global\liftamount=\z@
  \glet\maybeoffset=\relax}
\def\nudge#1{\ifinmeasureCD\else
  \global\advance\nudgeamount#1\relax
  \global\let\maybeoffset\dooffset\fi}
\def\lift#1{\ifinmeasureCD\else
  \global\advance\liftamount#1\relax
  \global\let\maybeoffset\dooffset\fi}

\def\findarrdp{\@rrdp=\CDvarrlen
  \ifnum\sp@ncnt@>1
    \advance\@rrdp by \CDstrutlen
    \multiply\@rrdp by \sp@ncnt@
    \advance\@rrdp by -\CDstrutlen \fi
 }

\def\varrow#1#2#3{\ifnum\sp@ncnt>\@ne 
     \sp@ncnt@=\sp@ncnt\relax\fi
  \findarrdp \@rrwd=\z@ 
  \kern\c@l
   \hbox to \z@{\options
   \@rrdp=\scalefactor\@rrdp
    \hss\vplace{#1}\llabel{#2}\rlabel{#3}\hss}%
  \global\advance\col by \@ne \setsp@n\@ne\@ne
  }

\def\novarrow{\varrow\vfill{}{}}

\def\tweenarrows#1{\findarrwd \findarrdp \setsp@n\@ne\@ne
  \rlap{\options\cplace{#1}}}

\def\usarrow #1#2#3{\dr@tn=\@ne
  \findarrwd \findarrdp \setsp@n\@ne\@ne 
  \rlap{\options\cplace{#1}\nwlabel{#2}\selabel{#3}}%
  \dr@tn=\z@}
\def\dsarrow #1#2#3{\dr@tn=\tw@
  \findarrwd \findarrdp \setsp@n\@ne\@ne 
  \rlap{\options\cplace{#1}\swlabel{#2}\nelabel{#3}}%
  \dr@tn=\z@}
 \def\@rrow#1{\csname #1@rrow\endcsname}
 \def\R@rrow{\harrow \rtarrfill}
 \def\L@rrow{\harrow \ltarrfill}
 \def\V@rrow{\varrow \dnarrfill}
 \def\A@rrow{\varrow \uparrfill}
 \def\SE@rrow{\dsarrow \SEarrow}
 \def\NW@rrow{\dsarrow \NWarrow}
 \def\SW@rrow{\usarrow \SWarrow}
 \def\NE@rrow{\usarrow \NEarrow}
 \def\DS@rrow{\dsarrow \dnslope}
 \def\US@rrow{\usarrow \upslope}
 \def\upslope{\find@xyargs
       \@linelen=\unitlength\@sline}
 \def\dnslope{\find@xyargs\@yarg=-\@yarg\relax
       \@linelen=\unitlength\@sline}

\newtoks\optionlist 
\optionlist={}
\let\options\relax
\def\dooptions{\the\optionlist\global\optionlist={}%
  \glet\options=\relax}
\def\option#1{\ifinmeasureCD\else
  \glet\options=\dooptions
  \global\optionlist=\expandafter{\the\optionlist\relax#1}\fi}
\def\wider#1{\ifinmeasureCD\else
   \option{\advance\@rrwd by #1}\fi}
\def\deeper#1{\ifinmeasureCD\else
   \option{\advance\@rrdp by #1}\fi}

{\def\\{\global\let\sptoken= }\\ }

\def\CR{\futurelet\nexttok\testCR}
\def\testCR{\ifx\nexttok\sptoken
   \let\next\eatspaceCR\else\let\next\CDCR\fi\next}
\def\eatspaceCR#1 {\CR}
\def\measureCR{\ifx\nexttok\endmeasure\let\nextCR\relax
    \else\let\nextCR\drop\fi\nextCR}
\def\setCR{\ifodd\row
  \ifx\nexttok\endsetCD\else\hend\betweenCDrows\vbegin\fi
  \else\vend\betweenCDrows\hbegin\fi}

\countdef\dimenc@unt=11
\def\CD#1\endCD{
   \begingroup\let\\=\CR
  \m@th\offinterlineskip
   \measure#1\endmeasure\null\,\vcenter{\setCD#1\endsetCD}\,
   \endgroup
    }

\ifx\@clnwd\undefined \nol@g\else\catcode`\ =14\relax\fi
 \font\@linefnt=line10 
 \newcount\@tempcnta
 \newcount\@tempcntb
 \newdimen\@tempdima
 \newdimen\@tempdimb
 \newdimen\@wholewidth
 \newdimen\@halfwidth
   \@wholewidth\fontdimen8\@linefnt \@halfwidth .5\@wholewidth
 \newdimen\unitlength
 \newcount\@xarg
 \newcount\@yarg
 \newcount\@yyarg
 \newbox\@linechar
 \newdimen\@linelen
 \newdimen\@clnwd
 \newdimen\@clnht
 \newif\if@negarg
 
 \def\@whilenoop#1{}

 \def\@whiledim#1\do #2{\ifdim #1\relax#2\@iwhiledim{#1\relax#2}\fi}

 \def\@iwhiledim#1{\ifdim #1\let\@nextwhile=\@iwhiledim 
         \else\let\@nextwhile=\@whilenoop\fi\@nextwhile{#1}}

 \def\@sline{\ifnum\@xarg< 0 \@negargtrue \@xarg -\@xarg \@yyarg -\@yarg
   \else \@negargfalse \@yyarg \@yarg \fi
 \ifnum \@yyarg >0 \@tempcnta\@yyarg \else \@tempcnta -\@yyarg \fi
 \ifnum\@tempcnta>6 \@badlinearg\@tempcnta0 \fi
 \ifnum\@xarg>6 \@badlinearg\@xarg 1 \fi
 \setbox\@linechar\hbox{\@linefnt\@getlinechar(\@xarg,\@yyarg)}%
 \ifnum \@yarg >0 \let\@upordown\raise \@clnht\z@
    \else\let\@upordown\lower \@clnht \ht\@linechar\fi
 \@clnwd=\wd\@linechar
 \if@negarg \hskip -\wd\@linechar \def\@tempa{\hskip -2\wd\@linechar}\else
      \let\@tempa\relax \fi
 \@whiledim \@clnwd <\@linelen \do
   {\@upordown\@clnht\copy\@linechar
    \@tempa
    \advance\@clnht \ht\@linechar
    \advance\@clnwd \wd\@linechar}%
 \advance\@clnht -\ht\@linechar
 \advance\@clnwd -\wd\@linechar
 \@tempdima\@linelen\advance\@tempdima -\@clnwd
 \@tempdimb\@tempdima\advance\@tempdimb -\wd\@linechar
 \if@negarg \hskip -\@tempdimb \else \hskip \@tempdimb \fi
 \multiply\@tempdima \@m
 \@tempcnta \@tempdima \@tempdima \wd\@linechar \divide\@tempcnta \@tempdima
 \@tempdima \ht\@linechar \multiply\@tempdima \@tempcnta
 \divide\@tempdima \@m
 \advance\@clnht \@tempdima
 \ifdim \@linelen <\wd\@linechar
    \hskip \wd\@linechar
   \else\@upordown\@clnht\copy\@linechar\fi}
 
 \def\@getlinechar(#1,#2){\@tempcnta#1\relax\multiply\@tempcnta 8
 \advance\@tempcnta -9 \ifnum #2>0 \advance\@tempcnta #2\relax\else
 \advance\@tempcnta -#2\relax\advance\@tempcnta 64 \fi
 \char\@tempcnta}
 
 \def\vector(#1,#2)#3{\@xarg #1\relax \@yarg #2\relax
 \@tempcnta \ifnum\@xarg<0 -\@xarg\else\@xarg\fi
 \ifnum\@tempcnta<5\relax
 \@linelen=#3\unitlength
 \ifnum\@xarg =0 \@vvector 
   \else \ifnum\@yarg =0 \@hvector \else \@svector\fi
 \fi
 \else\@badlinearg\fi}
 
 \def\@svector{\@sline
 \@tempcnta\@yarg \ifnum\@tempcnta <0 \@tempcnta=-\@tempcnta\fi
 \ifnum\@tempcnta <5
   \hskip -\wd\@linechar
   \@upordown\@clnht \hbox{\@linefnt  \if@negarg 
   \@getlarrow(\@xarg,\@yyarg) \else \@getrarrow(\@xarg,\@yyarg) \fi}%
 \else\@badlinearg\fi}
 
 \def\@getlarrow(#1,#2){\ifnum #2 =\z@ \@tempcnta='33\else
 \@tempcnta=#1\relax\multiply\@tempcnta \sixt@@n \advance\@tempcnta
 -9 \@tempcntb=#2\relax\multiply\@tempcntb \tw@
 \ifnum \@tempcntb >0 \advance\@tempcnta \@tempcntb\relax
 \else\advance\@tempcnta -\@tempcntb\advance\@tempcnta 64
 \fi\fi\char\@tempcnta}
 
 \def\@getrarrow(#1,#2){\@tempcntb=#2\relax
 \ifnum\@tempcntb < 0 \@tempcntb=-\@tempcntb\relax\fi
 \ifcase \@tempcntb\relax \@tempcnta='55 \or 
 \ifnum #1<3 \@tempcnta=#1\relax\multiply\@tempcnta
 24 \advance\@tempcnta -6 \else \ifnum #1=3 \@tempcnta=49
 \else\@tempcnta=58 \fi\fi\or 
 \ifnum #1<3 \@tempcnta=#1\relax\multiply\@tempcnta
 24 \advance\@tempcnta -3 \else \@tempcnta=51\fi\or 
 \@tempcnta=#1\relax\multiply\@tempcnta
 \sixt@@n \advance\@tempcnta -\tw@ \else
 \@tempcnta=#1\relax\multiply\@tempcnta
 \sixt@@n \advance\@tempcnta 7 \fi\ifnum #2<0 \advance\@tempcnta 64 \fi
 \char\@tempcnta}
\catcode`\ =10

\dol@g 
\catcode`\@=\active
\LaTeXfeatures

%% file: two_theorems_2405.bbl
\begin{thebibliography}{GKNP}

\bibitem[Ba]{Ba}
R.~Baer, {\it Engelsche Elemente Noetherscher Gruppen}, Math. Ann.
{\bf 133} (1957), 256--270.

\bibitem[BGGKPP]{BGGKPP}
 T.~Bandman, F.~Grunewald, G.-M.~Greuel, B.~Kunyavskii, G.~Pfister,
 E.~Plotkin,
{\it Two-variable Identities for Finite Solvable Groups},
C.R.~Acad. Sci. Paris, Ser. I, {\bf 337},  (2003), 581--586.


\bibitem[GS]{GS} E.~Golod, I.~Shafarevich, On the class field tower,
Izv. Akad. Nauk. SSSR Mat. Ser. 28 (1964), 261-272.

\bibitem [J]{J} N.~Jacobson, {\it Structure of rings},
American Math. Soc.Colloq. Publ., 37, 1964.

\bibitem [Ka]{Ka} L.~Kaluznin, {\it Uber gewisse Beziehungen
zwischen einer Gruppe und ihren Automorphismen}, Ber.
Mathematiker-Tagung, (1953) 164--172.

\bibitem [Ko]{Ko} A.I.~Kostrikin. Around Burnside, Springer-Verlag, Berlin- New York,
1990.

\bibitem [Ma]{Ma} A.I.~Malcev, {\it Faithful matrix
representations of infinite groups}, Math. Sbornic, {\bf 8}, (3),
(1940), 405--422.

\bibitem [Ne]{Ne} H.~Neumann, {\it Varieties of Groups},
Springer-Verlag, New York, 1967.

\bibitem[Pla]{Pla}
V.~P.~Platonov, {\it Engel elements and radical in PI-algebras and
topological groups}, DAN SSSR, {\bf 161}, (2), 289--291.

\bibitem[Pl1]{Pl1}
B.~I.~Plotkin, {\it Radical groups}, Mat.\ Sb.\ N.S.\ {\bf 37(79)}
(1955), 507--526; English transl.\ in Amer.\ Math.\ Soc.\ Transl.\
(2) {\bf 17} (1961), 9--28.

\bibitem[Pl2]{Pl2}
B.~I.~Plotkin, {\it Radical and nil-elements in groups}, Izvestija
Vuzov, ser. mat. (1958) (Russian).

\bibitem[Pl3]{Pl3}
B.~I.~Plotkin,{\it Groups of automorphisms of algebraic systems}.
Moscow, Nauka, 1966, Walter Nordhoff, Dordreht, 1971


\bibitem [Pr]{Pr} C.~Procesi, {\it The Burnside problem},
Journal of Algebra, {\bf 4}, (1966), 421--425.

\bibitem [Su]{Su} D. A.~Suprunenko, Matrix Groups, Amer. Math. Soc., Providence, 1976,
252 pp.

\bibitem[Th]{Th}
J.~Thompson, {\it Non-solvable finite groups all of whose local
subgroups are solvable}, Bull.\ Amer.\ Math.\ Soc.\ {\bf 74}
(1968), 383--437.

\bibitem [To1]{To1} A.~Tokarenko, {\it On radicals and stability
in linear groups over commutative rings  }, Sibirskii
Matematicheskii Zhurnal, (1968), vol. 9, (1), 165--176.


\bibitem [To2]{To2} A.~Tokarenko, {\it About linear groups over
rings}, Sibirskii Matematicheskii Zhurnal, (1968), vol. 9, (4),
951--959.



\bibitem[Wi]{Wi}
J.~S.~Wilson, {\it Two-generator conditions for residually finite
groups}, Bull.\ London Math.\ Soc.\ {\bf 23} (1991), 239--248.

\bibitem[WZ]{WZ}
J.~S.~Wilson and E.~Zelmanov, {\it Identities for Lie algebras of
pro-$p$ groups}, J.\ Pure Appl.\ Algebra {\bf 81} (1992),
103--109.

\end{thebibliography}
